# Why minimax is not that pessimistic

A. Fraysse*

**Abstract**

In nonparametric statistics an optimality criterion for estimation procedures is provided by the minimax rate of convergence. However this classical point of view is subject to controversy as it requires to look for the worst behaviour reached by an estimation procedure in a given space. The purpose of this paper is to show that this is not justified as the minimax risk often coincides with a generic one. We are here interested in the rate of convergence attained by some classical estimators on almost every, in the sense of prevalence, function in a Besov space.

## 1 Introduction

Since its introduction in the seventeen's, nonparametric estimation has taken a large place in the work of mathematical or signal processing communities. Often a signal has too many components, in the case for instance of densities or curve images, to allow classical studies upon finite dimensional spaces to give accurate results. But which kind of estimator can be more appropriate in nonparametric cases?

---

*LSS, CNRS UMR, University Paris-Sud, Supelec, 3 rue Joliot-Curie, 91190 Gif-Sur-Yvette, FRANCE. This work was performed when the author was at LTCI, Telecom ParisTech. Email: fraysse@lss.supelec.fr

This question raised a lot of definitions and discussions in statistical community. How can two estimators be compared when they point out infinite dimensional objects and what kind of optimal behaviour can be expected. One of the most common way to test the performance of a procedure is to compare its convergence rate with an optimal one given by minimax theory. Nonetheless, this technique comes from a particular definition which can be subject to controversy. The main drawback is the pessimist point of view of this theory, which looks for the worst rate of estimation obtained in a given space. In this paper, we introduce a new test of the risk, obtained thanks to genericity theories, and which shows that the minimax rate of estimation should not be as pessimistic as believed.

In the minimax paradigm one supposes that a function $f$ belongs to a certain space $\Theta$, which can be, for instance a Sobolev or Hölder space linked to some regularity properties and one defines a *risk*, or loss function thanks to a pseudo-distance on $\Theta$, denoted $R(.,.)$. Given a radius $C > 0$ and an estimation procedure $\hat{}_n$, depending on the model and of a data parameter $n$, the maximal risk of $\hat{}_n$ on $\Theta_C$ is then defined by:

$$R^n(\hat{}_n) = \sup_{f \in \Theta_C} \mathbb{E}(R(\hat{f}_n, f)), \qquad (1)$$

where $\Theta_C$ is a closed ball in $\Theta$ with radius $C > 0$.

If $\mathcal{T}_n$ denotes the set of all estimation procedures defined on $\Theta$ the minimax risk on $\Theta_C$ is given by :

$$R^n(\Theta) = \inf_{\hat{f} \in \mathcal{T}_n} \sup_{f \in \Theta_C} \mathbb{E}(R(\hat{f}_n, f)). \qquad (2)$$

This minimax risk gives an optimal bound over the function class $\Theta_C$. It is thus natural for estimation procedures to attempt to reach this risk, at least asymptotically when $n$ tends to infinity.

The main drawback of the minimax theory is that we are looking for the maximum risk on a function space, thus for the worst behaviour. This point of view seems pessimistic and can be not generic enough as it can be used to merge estimation procedures. Indeed, the worst case could be a misleading one and a method can be rejected although it is a good one for a lot of functions. The purpose of this paper is to show that it is not the case and that in fact minimax rate corresponds to a generic one. The purpose of this paper is thus twofold. We introduce a new test of estimation performances, based on generic properties. Thanks to this definition, we show that the minimax risk coincide with an "almost every" one.

Let us first introduce what is meant by almost every function. In a finite dimensional space, we say that a property holds almost everywhere if the set of points where it is not true is of vanishing Lebesgue measure. The Lebesgue measure has here a preponderant role, as it is the only $\sigma$-finite and translation invariant measure. Unfortunately, no measure shares those properties in infinite dimensional Banach spaces. A way to recover a natural "almost every" notion in infinite vector spaces is thus defined as follows by J. Christensen in 1972 see [2, 4, 12].

**Definition 1.** *Let $V$ be a complete metric vector space. A Borel set $A \subset V$ is Haar-null (or shy) if there exists a compactly supported probability measure $\mu$ such that*

$$\forall x \in V, \quad \mu(x + A) = 0. \tag{3}$$

*If this property holds, the measure $\mu$ is said to be transverse to $A$.*

*A subset of $V$ is called Haar-null if it is contained in a Haar-null Borel set. The complement of a Haar-null set is called a prevalent set.*

As it can be seen in the definition of prevalence, the main issue in proofs is to construct transverse measures to a Borel Haar-null set. We remind here two classical ways to construct such measures.

**Remarks.**  1. A finite dimensional subspace of $V$, $P$, is called a probe for a prevalent set $T \subset V$ if the Lebesgue measure on $P$ is transverse to the complement of $T$.
   Those measures are not compactly supported probability measures. However one immediately checks that this notion can also be defined the same way but stated with the Lebesgue measure defined on the unit ball of $P$. Note that in this case, the support of the measure is included in the unit ball of a finite dimensional subspace. The compactness assumption is therefore fulfilled.

2. If $V$ is a function space, a probability measure on $V$ can be defined by a random process $X_t$ whose sample paths are almost surely in $V$. The condition $\mu(f + A) = 0$ means that the event $X_t - f \in A$ has probability zero. Therefore, a way to check that a property $\mathcal{P}$ holds only on a Haar-null set is to exhibit a random process $X_t$ whose sample paths are in $V$ and is such that

$$\forall f \in V, \ a.s. \quad X_t + f \ \text{does not satisfy } \mathcal{P}.$$

The following results enumerate important properties of prevalence and show that these notions supply a natural generalization of "zero measure" and "almost every" in finite-dimensional spaces, see [2, 4, 12].

**Proposition 1.**
- *If $S$ is Haar-null, then $\forall x \in V$, $x + S$ is Haar-null.*
- *If $\dim(V) < \infty$, $S$ is Haar-null if and only if $\operatorname{meas}(S) = 0$ (where meas denotes the Lebesgue measure).*
- *Prevalent sets are dense.*
- *The intersection of a countable collection of prevalent sets is prevalent.*
- *If $\dim(V) = \infty$, compact subsets of $V$ are Haar-null.*

As we can see from the properties of prevalent sets, this theory provides a natural generalization of the finite dimensional notion of almost every. Since its definition, this theory has been mainly used in the context of differential geometry [12] and regularity type properties [11]. A classical example is given in [11], where it is proved that the set of nowhere differentiable functions is prevalent in the space of continuous functions.

Using this theory, a natural way to exhibit a risk for an estimating procedure is to look at the risk reached on almost every function of $\Theta$, in the sense of prevalence.

As the minimax theory has been widely studied, a large class of results exist in different function spaces and with different losses. Historically, the first one is the result of Pinsker [21] which shows that suitable linear estimators reach the optimal $L^2$ risk rate on $L^2$ Sobolev classes. If the risk function is given by an $L^p$ norm, [13, 3] show that, under certain conditions, kernel estimators are optimal in the sense of minimax theory in the same function spaces. More recent results, such as those of [20], stated that linear estimators cannot reach the optimal bound in nonlinear regression, as soon as we take the $L^p$ risk and Sobolev classes.

In this paper we focus on Besov spaces and take the general $L^p$ loss function. The interest of studying Besov spaces is motivated by its practical use in approximation theory and its theoretical simplicity in terms of wavelet expansions. Furthermore, in the theoretical point of view, they also generalize some classical function space, such as Hölder and $L^2$ Sobolev spaces.

In those Besov spaces we study the performances in terms of generic approximation of two classical estimation procedures in both white noise model and density estimation problem. A second result gives the generic rate of estimation for larger families of procedure in the case of the white noise model.

## 2 Models and estimation procedures

In the following, we consider two classical estimation problems. The first one is given by the Gaussian white noise model. Following the definition of [13], we suppose that we observe $Y_t$ such that

$$dY_t = f(t)dt + \frac{1}{\sqrt{n}}dW_t, \quad t \in (0,1)^d, \qquad (4)$$

where $dW_t$ stands for the $d$-dimensional Wiener measure, $n$ is known and $f$ is the unknown function to be estimated.

The second theoretical framework in which our theory can be applied is the problem of density estimation. In this case, the model is given by a sequence $X_1, \ldots, X_n$ of independent and identically distributed random variables of unknown density $f$ on $\mathbb{R}$.

The estimation procedures that we deal with are defined thanks to a base decomposition of the function to be estimated. To define them, we first introduce the wavelet decomposition. In our framework, those bases allow both to define function spaces and estimation procedures. It provide thus a key tool to introduce our results.

The wavelet transform is a powerful approximation tool largely used in statistics applied to signal processing, thanks to its properties of localization in time and frequency domains. Indeed, this property allows to reconstruct a signal with few coefficients. Its use in statistical communities and the development of wavelet based estimators are thus natural, as introduced in [18].

To define wavelets, we refer to [6] where it is proved that for $r$ large enough there exists $2^d - 1$ functions $\psi^{(i)}$ with compact support and which belong to $C^r$. Furthermore each $\psi^{(i)}$ has $r$ vanishing moments and the set of functions $\{\psi_{j,k}^{(i)} = 2^{dj/2}\psi^{(i)}(2^j x - k), \quad j \in \mathbb{Z}, \ k \in \mathbb{Z}^d, i \in \{1, ..., 2^d - 1\}\}$ forms an orthonormal basis of $L^2(\mathbb{R}^d)$. It is also noticed in [19] that wavelets provide unconditional bases of $L^p(\mathbb{R}^d)$ as far as $1 < p < \infty$. Taking periodized wavelets allow to restraint our properties to $[0,1]^d$.

Thus any function $f \in L^p$ can be written as

$$f = \sum_{i,j,k} c_{j,k}^{(i)} \psi_{j,k}(x) \tag{5}$$

where

$$c_{j,k}^{(i)} = 2^{jd/2} \int f(x) \psi^{(i)}(2^j x - k) dx. \tag{6}$$

We can notice that we stand in isotropic cases. Thus the direction of the wavelets is not involved and in the following, we omit the directional index $i$.

As the collection of $\{2^{dj/2}\psi^{(i)}(2^j x - k), \quad j \in \mathbb{N}, \ k \in \{0, \ldots, 2^j - 1\}^d, i = 1, ..., 2^d - 1\}$ form an orthonormal basis of $L^2([0,1]^d)$, observing the whole trajectory of $Y_t$ in (4) is equivalent to treat the following problem, in which is observed $(y_{j,k})_{j \in \mathbb{N}, k \in \{0,...,2^j-1\}^d} \in \ell^2(\mathbb{N}^{d+1})$ such that $\forall j, k$,

$$y_{j,k} = \theta_{j,k} + \frac{1}{\sqrt{n}} v_{j,k}, \tag{7}$$

where $y_i = \int \psi_{j,k} dY$, $v_i$ are i.i.d. Gaussian random variables and $(\theta_{j,k})$ is the sequence to be estimated.

In terms of density estimation, one can also notice that the density function to be estimated $f$ can be represented in terms of wavelets $f = \sum \beta_{j,k} \psi_{j,k}$. In this case, the purpose is to find a sequence $(\hat{\beta}_{j,k})_{j,k}$ approximating $(\beta_{j,k})_{j,k}$.

Furthermore wavelets are useful as they provide a simple characterization of Besov spaces.

Homogeneous Besov spaces are characterized, for $p, q > 0$ and $s \in \mathbb{R}$, by

$$f \in B_p^{s,q}([0,1]^d) \iff \exists C > 0 \quad \sum_j \left( \sum_{j \geq 0, k \in \{0, \ldots, 2^j - 1\}^d} |c_{j,k}|^p 2^{(sp - d + \frac{pd}{2})j} \right)^{q/p} \leq C. \tag{8}$$

This characterization is independent from the chosen wavelet has $r$ vanishing moments, with $r \geq s$.

We also denote by $B_{p,c}^{s,q}(\mathbb{R}^d)$ the closed ball in $B_p^{s,q}(\mathbb{R}^d)$ of radius $c > 0$.

In terms of wavelets approximation, or in any base, the most natural and classical way to define estimators reachable is given by linear estimation.

**Definition 2.** *Suppose that we stand in the model (7). Linear estimators $\hat{f}_n^L$ are constructed by*

$$\hat{f}_n^L(x) = \sum_{j,k} \hat{\theta}_{j,k}^{(n)} \psi_{j,k}(x), \tag{9}$$

*where*

$$\hat{\theta}_{j,k}^{(n)} = \lambda_{j,k}^{(n)} y_{j,k}. \tag{10}$$

Parameters $(\lambda_{j,k}^{(n)})_{j,k}$ can be seen as smoothing weights depending of the problem. Those weights can be of different natures. Classical ones are:

- *Projection weights:* $\lambda_{j,k}^{(n)} = \mathbb{1}_{j < m_n}$.
- *Pinsker weights:* $\lambda_{j,k}^{(n)} = (1 - (\frac{j}{m_n})^\alpha)_+$,

where $m_n$ is an increasing function of $n$.

**Definition 3.** *Suppose that we stand in the model of density estimation. In this case, a linear estimator of the density $f$ is constructed by taking*

$$\hat{\beta}_{j,k} = \frac{1}{n} \sum_{i=1}^{n} <\psi_{j,k}, X_i>. \tag{11}$$

*And*

$$\hat{f}_n^L = \sum_{j,k} \hat{\beta}_{j,k} \psi_{j,k}. \tag{12}$$

The localization property of wavelet expansions is such that a given signal may have a sparse representation in those bases. Thus a natural estimation procedure in the white noise model, defined in [7] and ever since widely used in the signal community is to take away small wavelet coefficients. This is the principle of wavelet thresholding.

**Definition 4.** *Suppose that we stand in the case of white noise model (7). The wavelet thresholding procedure is then defined by*

$$\hat{f}_n^T(x) = \sum_{j=0}^{j(n)} \sum_k \beta_{j,k}^T \psi_{j,k}. \tag{13}$$

*Here the weights are given by:*

$$\beta_{j,k}^T = y_{j,k} \mathbb{1}_{\{|y_{j,k}| \geq \kappa t_n\}}, \tag{14}$$

*in the case of hard thresholding, or*

$$\beta_{j,k}^T = \text{sign}(y_{j,k})(|y_{j,k}| - \kappa t_n)_+, \tag{15}$$

*for the soft thresholding. Here,*

$$t_n = \sqrt{\frac{\log n}{n}}, \tag{16}$$

*stands for the universal thresholding and $j(n)$ is such that*

$$2^{-j(n)} \leq \frac{\log n}{n} < 2^{-j(n)+1},$$

*$\kappa$ being a constant large enough.*

Once again, in the model of density estimation wavelet thresholding is obtained thanks to a slight modification of the previous definition.

**Definition 5.** *Suppose that we stand in the problem of density estimation, and let $\hat{\beta}_{j,k}$ be the coefficients defined in (11). Thus the density estimator by wavelet thresholding is given by*

$$\hat{f}_n^T = \sum_{j=1}^{j(n)} \sum_k \hat{\beta}_{j,k} \mathbb{1}_{\{|\hat{\beta}_{j,k}| > t_n\}}. \tag{17}$$

*Where*

$$t_n = \sqrt{\frac{\log n}{n}}, \quad (18)$$

*is the universal thresholding and $j(n)$ is such that*

$$2^{-j(n)} \leq \frac{\log n}{n} < 2^{-j(n)+1}$$

Those estimators all belong to larger classes of estimation procedures, namely the classes of limited and elitist rules.

**Definition 6.** *Suppose that we stand in the white noise model (7). Let us consider the class of shrinkage estimators*

$$\mathcal{F}_n = \left\{ \hat{f}_n = \sum \gamma_{j,k} y_{j,k} \psi_{j,k};\ \gamma_{j,k} \in [0,1],\ measurable \right\}.$$

*Thus we say that $\hat{f}_n \in \mathcal{F}_n$ is a limited rule if there exist a deterministic function $\lambda_n$ and a constant $a > 0$ such that for every $j, k$:*

$$\gamma_{j,k} > a \Rightarrow 2^{-j} > \lambda_n. \quad (19)$$

*In this case, one say that $\hat{f}_n$ belongs to the class $\mathcal{L}(\lambda_n, a)$.*

*We say that $\hat{f}_n \in \mathcal{F}_n$ is an elitist rule if there exist a deterministic function $\lambda_n$ and a constant $a > 0$ such that for every $j, k$:*

$$\gamma_{j,k} > a \Rightarrow |y_{j,k}| > \lambda_n. \quad (20)$$

*Thus, $\hat{f}_n$ belongs to the class $\mathcal{E}(\lambda_n, a)$.*

One can easily see that linear estimators introduced in Definition 2 are limited rules and that thresholding algorithms, hard and soft thresholding, or some Bayesian procedures with Gaussian prior are elitist rules.

# 3 Statement of main results

Let us recall minimax results in Besov spaces. Taking the $L^p$ norm, where $1 \leq p < \infty$, as the loss function, we know from [9], for the white noise model or in the case of density estimation, that the minimax lower bound in closed balls in Besov spaces in given by the following proposition.

**Proposition 2.** *Let $1 \leq r \leq \infty$, $1 \leq p < \infty$ and $s > \frac{d}{r}$. Then, there exists $C > 0$ such that*

$$R^n(B_{r,c}^{s,\infty}) = \inf_{T_n} \sup_{f \in B_{r,c}^{s,\infty}} \mathbb{E}\|T_n - f\|_{L^p}^p \geq C r_n(s,r,p) \qquad (21)$$

*where*

$$r_n(s,r,p) = \begin{cases} n^{-\frac{ps}{2s+d}} & \text{if } r > \frac{dp}{2s+d}, \\ \left(\frac{n}{\log n}\right)^{-\frac{p(s-\frac{d}{r}+\frac{d}{p})}{2(s-\frac{d}{r})+d}} & \text{else.} \end{cases} \qquad (22)$$

Let us now check what is known concerning estimation procedures that we deal with. Although it is proved in [9] that thresholding procedures reach asymptotically the optimal rate up to a logarithmic correction, it is not always the case for linear procedures. As it can be seen in [8], with $L^2$ risk, linear estimators do not attain the minimax rate when studied functions have a sparse representation in a given base. In fact, the following proposition is proved in [8] which gives the optimal rate that can be reached in this case.

**Proposition 3.** *Let $1 \leq r \leq \infty$, $1 \leq p < \infty$ and $s > \frac{d}{r}$. There exist $C > 0$ such that*

$$R_{lin}^n(B_{r,c}^{s,\infty}) = \inf_{\tilde{T}_n \, linear} \sup_{f \in B_{r,c}^{s,\infty}} \mathbb{E}\|\tilde{T}_n - f\|_{L^p}^p \geq C \tilde{r}_n(s,r,p) \qquad (23)$$

*where*

$$\tilde{r}_n(s,r,p) = \begin{cases} n^{-\frac{ps}{2s+d}} & \text{if } r > p \\ \left(\frac{n}{\log n}\right)^{-\frac{ps'}{2s'+d}} & \text{else,} \end{cases} \qquad (24)$$

*and $s' = s - \frac{d}{r} + \frac{d}{p}$.*

We see in the following theorem that Proposition 3 remains true if we replace the risk maximum by the risk reached for almost every function. We also prove that in the same context Proposition 2 is true for thresholding algorithms up to a logarithmic term. We say in the following that $a_n \approx b_n$ if $\frac{\log a_n}{\log b_n} \to 1$.

**Theorem 1.** *Let $1 \leq r \leq \infty$, $1 \leq p < \infty$ and $s > \frac{d}{r}$. Then, in the context of (4) or for the problem of density estimation:*

- *For a suitable linear estimator $\hat{f}_n^L$ as in Definition 2 and for almost every function $f$ in $B_r^{s,\infty}([0,1]^d)$,*

$$\mathbb{E}\|\hat{f}_n^L - f\|_{L^p}^p \approx n^{-\alpha p}, \qquad (25)$$

where

$$\alpha = \begin{cases} \frac{s}{2s+d} & \text{if } r \geq p, \\ \frac{s-\frac{d}{r}+\frac{d}{p}}{2(s-\frac{d}{r}+\frac{d}{p})+d} & \text{else.} \end{cases} \quad (26)$$

- For almost every function in $B_r^{s,\infty}([0,1]^d)$, and for thresholding estimator $\hat{f}_n^T$

$$\mathbb{E}\|\hat{f}_n^T - f\|_{L^p}^p \approx \left(\frac{n}{\log n}\right)^{-\alpha p} \quad (27)$$

where

$$\alpha = \begin{cases} \frac{s}{2s+d} & \text{if } r > \frac{pd}{2s+d} \\ \frac{(s-\frac{d}{r}+\frac{d}{p})}{2(s-\frac{d}{r})+d} & \text{else.} \end{cases} \quad (28)$$

This theorem can be extended in term of shrinkage classes thanks to the following result.

**Theorem 2.** *Let $1 \leq r \leq \infty$, $1 \leq p < \infty$ and $s > \frac{d}{r}$. Then, in the context of (4) one has*

1. *For limited rules, and for almost every function $f$ in $B_r^{s,\infty}([0,1]^d)$, for every $C > 0$,*

$$\mathbb{E}\|\hat{f}_n - f\|_{L^p}^p \geq Cn^{-\alpha p}, \quad (29)$$

where

$$\alpha = \begin{cases} \frac{s}{2s+d} & \text{if } r \geq p, \\ \frac{s-\frac{d}{r}+\frac{d}{p}}{2(s-\frac{d}{r}+\frac{d}{p})+d} & \text{else.} \end{cases} \quad (30)$$

2. *For almost every function in $B_r^{s,\infty}([0,1]^d)$, and for elitist rules, for every $C > 0$,*

$$\mathbb{E}\|\hat{f}_n - f\|_{L^p}^p \geq C\left(\frac{n}{\log n}\right)^{-\alpha p} \quad (31)$$

where

$$\alpha = \begin{cases} \frac{s}{2s+d} & \text{if } r > \frac{pd}{2s+d} \\ \frac{(s-\frac{d}{r}+\frac{d}{p})}{2(s-\frac{d}{r})+d} & \text{else.} \end{cases} \quad (32)$$

The previous theorems are stated in terms of Besov spaces and wavelet based estimation. Nevertheless, they can be easily extended to other function spaces, such as Sobolev spaces, thanks to an adapted basis. One can also notice that all of our results are given in terms of polynomial rate of estimation.

# 4 Proofs of Theorem 1 and 2

Indeed, the proof of this theorem is based on the maximal space where an estimation procedure attains a given rate of convergence. For the sake of completeness, let us recall some basic facts upon the corresponding theory.

## 4.1 Maxiset theory

The maxiset theory introduced recently in [5, 16, 17] is an alternative way to compare different estimation procedures. In our case, it provides a crucial key to prove Theorem 1. The main idea is to look for the maximal space on which an estimator will reach a given rate, instead of searching an optimal rate for a given space.

**Definition 7.** *Let $\rho$ be a risk function and $(v_n)_{n \in \mathbb{N}}$ a sequence such that $v_n \to 0$. For $\hat{}_n$ a procedure, the maximal space associated to $\rho$, $v_n$ and a constant $T$ is given by*

$$MS(\hat{}_n, \rho, v_n, T) = \left\{ f;\ \sup_n v_n^{-1} \mathbb{E}(\rho(\hat{}_n, f)) < T \right\}. \tag{33}$$

Several improvements were made in nonparametric theory thanks to this idea. For instance, it is shown in [5] that, for the density estimation model the thresholding procedure is more efficient than the linear procedure, whose maxiset is given in [15]. And in the heteroscedastic white noise model, [22, 23] shown that thresholding procedures are better than linear estimators and as good as Bayesian procedures. In the case of white noise model, we recall the following result which is a particular case of [22].

**Proposition 4.** *Let $1 \leq p < \infty$, $1 \leq r < \infty$, $s > \frac{d}{r}$ and $\alpha \in (0,1)$. Let $\hat{f}_n^L$ be the estimator given in Definition 2.*
*For every $f$, we have the following equivalence:*
*There exists $c > 0$ such that for every $n \in \mathbb{N}$,*

$$\mathbb{E}\|\hat{f}_n^L - f\|_p^p \leq c m_n^{-\alpha p} \tag{34}$$

*if and only if $f \in B_p^{\alpha, \infty}$.*

Before stating the result associated with thresholding algorithms, we define new function spaces closely related to approximation theory. Those spaces, weak Besov spaces, defined in [5] are subsets of Lorentz spaces, and constitute a larger class than Besov spaces.

**Definition 8.** *Let $0 < r < p < \infty$. We say that a function $f = \sum_{j,k} c_{j,k} \psi_{j,k}$ belongs to $W(r,p)$ if and only if*

$$\sup_{\lambda > 0} \lambda^r \sum 2^{j(\frac{dp}{2}-d)} \sum_k \mathbb{1}_{\{|c_{j,k}|>\lambda\}} < \infty. \tag{35}$$

A fast calculation shows that the space $W(r,p)$ contains Besov spaces $B_r^{\beta,\infty}$ as soon as $\beta \geq \frac{d}{2}(\frac{p}{r} - 1)$.

The maxiset associated with the thresholding estimation procedure is given by a weak Besov space as proved in [5], and developed further in the heteroscedastic regression case in [16].

**Proposition 5.** *Let $1 \leq p < \infty$, $1 \leq r < \infty$, $s > \frac{d}{r}$ and $\tilde{\alpha} \in (0,1)$. Let $\hat{f}_n^T$ be the estimator defined by (4) and (15). Then for every $f$ we have the following equivalence:*
*$\exists K > 0$ such that $\forall n > 0$,*

$$\mathbb{E}\|\hat{f}_n^T - f\|_p^p \leq K \left(\sqrt{n \log(n)^{-1}}\right)^{-\tilde{\alpha}p} \tag{36}$$

*if and only if $f \in B_p^{\tilde{\alpha}/2,\infty} \cap W((1-\tilde{\alpha})p, p)$.*

Concerning shrinkage procedures, the following result from [1] also gives maxiset results for a large class of procedures.

**Proposition 6.** *Let $1 \leq p \leq \infty$ and $\sigma > 0$ be fixed.*

- *Let $\hat{f}_n$ be a limited rule in $\mathcal{L}(\lambda_n, a)$, with $a \in [0, 1[$ and $\lambda_n$ a non decreasing continuous function. Thus*

$$MS(\hat{f}_n, \|.\|_p^p, \lambda_n^{p\sigma}) \subset B_{p,\infty}^{\sigma}. \tag{37}$$

- *Let $\hat{f}_n$ be an elitist rule in $\mathcal{E}(\lambda_n, a)$, with $a \in [0, 1[$ and $\lambda_n$ a non decreasing continuous function. Thus*

$$MS(\hat{f}_n, \|.\|_p^p, \lambda_n^{p-\sigma}) \subset W(\sigma, p). \tag{38}$$

Furthermore, another important key result involving Besov spaces is the following proposition from [10].

**Proposition 7.** *Let us define the scaling function of a distribution $f$ by*

$$\forall p > 0 \quad s_f(p) = \sup\{s : f \in B_p^{s,\infty}\}. \tag{39}$$

Let $s_0$ and $p_0$ be fixed such that $s_0 - \frac{d}{p_0} > 0$. Outside a Haar-null set in $B_{p_0}^{s_0,\infty}(\mathbb{R}^d)$, we have:

$$s_f(p) = \begin{cases} s_0 & \text{if } p \leq p_0 \\ \frac{d}{p} + s_0 - \frac{d}{p_0} & \text{if } p \geq p_0. \end{cases} \quad (40)$$

One can check that a lower bound of this scaling function is given by Besov embeddings and interpolation theory, which can be found in [24]. This result thus state that we cannot have a better regularity than the one given by those embeddings. In our case, we will exploit this result by comparing those critical spaces with the maxiset associated to each procedure.

## 4.2 Generic risk for linear estimators

Let $1 \leq p < \infty$, $1 \leq r < \infty$ and $s > \frac{d}{r}$ be fixed. Denote

$$s' = s - \left(\frac{d}{r} - \frac{d}{p}\right)_+,$$

and

$$\alpha(s') = \frac{s'}{2s' + d}. \quad (41)$$

In this section, we prove the first part of Theorem 1. We define the linear estimator as in Definitions 2 and 3. A bias-variance compromise shows that for $r \geq p$, we have to take $m_n = n^{\frac{1}{2s+d}}$ whereas when $p > r$, the bias and the variance are compensated when $m_j = n^{\frac{1}{2(s-\frac{d}{r}+\frac{d}{p})+d}}$.

The upper bound is straightforward. Actually, from [7] we know that there exists $c > 0$ such that for any $n \in \mathbb{N}$, and for any $f \in B_r^{s,\infty}([0,1]^d)$,

$$\mathbb{E}\|\hat{f}_n^L - f\|_p^p \leq cn^{-\alpha(s')p}. \quad (42)$$

Let us now check the lower bound. In a first time, we have to show that for every $\varepsilon > 0$ fixed, the set

$$M(\varepsilon) = \left\{f \in B_r^{s,\infty}([0,1]^d);\ \exists c > 0\ \forall n \in \mathbb{N},\ \mathbb{E}(\|\hat{f}_n^L - f\|_{L^p}^p) < cn^{-(\alpha(s')+\varepsilon)p}\right\} \quad (43)$$

is a Haar null set.

Furthermore, with the particular form of $m_n$ taken, it coincides with the following set:

$$\left\{ f \in B_r^{s,\infty}([0,1]^d);\ \exists c > 0\ \forall n \in \mathbb{N},\ \mathbb{E}(\|\hat{f}_n^L - f\|_{L^p}^p) < c m_n^{-(s'+\varepsilon)p} \right\}.$$

But, by applying Proposition 4 we see that $M(\varepsilon)$ is included in $B_p^{s'+\varepsilon,\infty}([0,1]^d)$. And from Proposition 7, we know that this set is a Haar null Borel set of $B_r^{s,\infty}([0,1]^d)$.

We thus obtain that $\forall \varepsilon > 0$, the set

$$\left\{ f \in B_r^{s,\infty}([0,1]^d);\ \exists c > 0\ \forall n \in \mathbb{N},\ \mathbb{E}(\|\hat{f}_n^L - f\|_{L^p}^p) < c n^{-(\alpha(s')+\varepsilon)p} \right\}$$

is a Haar null set. This set can also be written,

$$\left\{ f \in B_r^{s,\infty}([0,1]^d);\ \limsup_{n \to \infty} \frac{\log(\mathbb{E}(\|\hat{f}_n^L - f\|_{L^p}^p))}{-p \log n} > \alpha(s') + \varepsilon \right\}.$$

Taking the countable union of those sets over a decreasing sequence $\varepsilon_n \to 0$, and the complementary we obtain that for almost every function in $B_r^{s,\infty}([0,1]^d)$,

$$\liminf_{n \to \infty} \frac{\log(\mathbb{E}(\|\hat{f}_n^L - f\|_{L^p}^p))}{-p \log n} \leq \alpha(s'). \tag{44}$$

Which induces the expected result.

## 4.3 Thresholding algorithms

In this part, we take the estimation procedures $\hat{f}_n^T$ given in Definition 4 and Definition 5.

Let us turn out our attention to the minimax rate of convergence. For this purpose, we write in the following

$$\tilde{\alpha}(s) = \begin{cases} \frac{2s}{2s+d} & \text{if } r > \frac{pd}{2s+d} \\ \frac{2(s-\frac{d}{r}+\frac{d}{p})}{2(s-\frac{d}{r})+d} & \text{else.} \end{cases} \tag{45}$$

The proof of the second point of Theorem 1 follows the same scheme as the previous one. In this case, the upper bound is given in

[9]. Thus we know that for every function in $B_r^{s,\infty}([0,1]^d)$, and for all $1 < p < \infty$,
$$\mathbb{E}(\|\hat{f}_n^T - f\|_{L^p}^p) < c\sqrt{\frac{n}{\log n}}^{-\tilde{\alpha}(s)p}.$$

In order to prove the lower bound, we use Proposition 5.

For every values of $\tilde{\alpha}$, let $0 < \varepsilon < 1 - \tilde{\alpha}$ be fixed, and $M(\varepsilon)$ be the set defined by
$$M(\varepsilon) = \left\{ f \in B_r^{s,\infty}([0,1]^d);\ \exists c > 0\ \forall n \in \mathbb{N},\ \mathbb{E}(\|\hat{f}_n^T - f\|_{L^p}^p) < c\sqrt{\frac{n}{\log n}}^{-(\tilde{\alpha}(s)+\varepsilon)p} \right\}.$$

Thanks to Proposition 5, this set $M(\varepsilon)$ is embedded in $B_p^{\frac{\tilde{\alpha}+\varepsilon}{2},\infty} \cap W((1-\tilde{\alpha}-\varepsilon)p, p)$.

The end of the proof is based on the following proposition.

**Proposition 8.** *Let us define the weak scaling function of a distribution $f$ by*
$$\forall p > 0 \quad \tilde{s}_f(p) = \sup\{\alpha :\ f \in W((1-\alpha)p, p)\}. \tag{46}$$

*Let $s$ and $r$ be fixed such that $s - \frac{d}{r} > 0$. Outside a Haar-null set in $B_r^{s,\infty}(\mathbb{R}^d)$, we have:*
$$\tilde{s}_f(p) = \begin{cases} \frac{2s}{2s+d} & \text{if } r > \frac{pd}{2s+d} \\ \frac{2(s-\frac{d}{r}+\frac{d}{p})}{2(s-\frac{d}{r})+d} & \text{else.} \end{cases} \tag{47}$$

*Proof*: In order to prove Proposition 8, let us prove that $W((1-\tilde{\alpha}-\varepsilon)p, p)$ is a Haar null Borel set in $B_r^{s,\infty}([0,1]^d)$. For this purpose, we define our transverse measure as the probe generated by the function $g$ defined by its wavelet coefficients:
$$d_{j,k} = \frac{2^{-(s-\frac{d}{r}+\frac{d}{2})j} 2^{-\frac{d}{r}J}}{j^a} \tag{48}$$
where $a = 1 + \frac{3}{r}$ and $0 \leq J \leq j$ and $K \in \{0, \ldots, 2^J - 1\}^d$ are such that
$$\frac{K}{2^J} = \frac{k}{2^j} \tag{49}$$

is an irreducible fraction. As it can be seen in [14], this function $g$ belongs to $B_r^{s,\infty}([0,1]^d)$. Let $f \in B_r^{s,\infty}([0,1]^d)$ be an arbitrary function and consider the affine subset

$$M = \{\alpha \in \mathbb{R} \ f + \alpha g \in W((1 - \tilde{\alpha} - \varepsilon)p, p)\}. \tag{50}$$

Suppose that there exist two points $\alpha_1$ and $\alpha_2$ in $M$. Thus $f + \alpha_1 g - (f + \alpha_2 g)$ belongs to $W((1 - \tilde{\alpha} - \varepsilon)p, p)$, and therefore there exists $c > 0$ such that

$$\|f + \alpha_1 g - (f + \alpha_2 g)\|_{W((1-\tilde{\alpha}-\varepsilon)p,p)} = \|(\alpha_1 - \alpha_2)g\|_{W((1-\tilde{\alpha}-\varepsilon)p,p)} \leq c. \tag{51}$$

As a fast calculation shows that

$$\forall \alpha > 0, \ \|\alpha g\|_{W(r,p)} = \alpha^r \|g\|_{W(r,p)} \tag{52}$$

we just have now to determine $\|g\|_{W(r,p)}$. Thanks to equation (35), this is equivalent to determine for every $t > 0$, the value of

$$2^{-(1-\tilde{\alpha}-\varepsilon)pt} \sum_{j \geq 0} 2^{j(\frac{dp}{2}-d)} \sum_k \mathbb{1}_{\{d_{j,k} > 2^{-t}\}}$$

But by definition of $g$, we have,

$$\frac{2^{-(s-\frac{d}{r}+\frac{d}{2})j} 2^{-\frac{d}{r}J}}{j^a} > 2^{-t} \Rightarrow (s - \frac{d}{r} + \frac{d}{2})j + \frac{d}{r}J \leq t,$$

which implies that

$$J \leq \frac{r}{d}t - (s - \frac{d}{r} + \frac{d}{2})\frac{r}{d}j. \tag{53}$$

Note that the condition $J \geq 0$ implies also that

$$j(s - \frac{d}{r} + \frac{d}{2}) \leq t. \tag{54}$$

We denote by $\tilde{t} = \frac{t}{s-\frac{d}{r}+\frac{d}{2}}$ and by $\tilde{\tilde{t}} = \frac{t}{s+\frac{d}{2}}$. Thus we have, for every $t > 0$,

$$\|g\|_{W((1-\tilde{\alpha}-\varepsilon)p,p)} \geq 2^{-(1-\tilde{\alpha}-\varepsilon)pt} \sup_{0 \leq j \leq \tilde{t}} 2^{j(\frac{dp}{2}-d)} \sum_{J=0}^{j \wedge [\frac{r}{d}t - (s-\frac{d}{r}+\frac{d}{2})\frac{r}{d}j]} 2^{dJ}$$

$$\geq 2^{-(1-\tilde{\alpha}-\varepsilon)pt} \sup \left( \sup_{0 \leq j \leq \frac{t}{s+\frac{d}{2}}} 2^{j(\frac{dp}{2}-d)} \sum_{J=0}^{j} 2^{dJ}, \sup_{\frac{t}{s+\frac{d}{2}}+1 \leq j \leq \tilde{t}} 2^{j(\frac{dp}{2}-d)} \sum_{J=0}^{[\frac{r}{d}t-(s-\frac{d}{r}+\frac{d}{2})\frac{r}{d}j]} 2^{dJ} \right)$$

$$\geq \frac{2^{-(1-\tilde{\alpha}-\varepsilon)pt}}{2^d - 1} \sup \left( \sup_{0 \leq j \leq \tilde{\tilde{t}}} 2^{\frac{dpj}{2}}(1 - 2^{-jd}), \sup_{\tilde{\tilde{t}} < j \leq \tilde{t}} 2^{j(\frac{dp}{2}-d)}(2^{rt}2^{-jr(s+\frac{d}{2}-\frac{d}{r})} - 1) \right)$$

Merging this result with (51) together with (52), we obtain that, if there exist $\alpha_1$ and $\alpha_2$ in $M$ then they satisfy that for every $t \geq 0$ and $0 \leq j \leq \tilde{t}$,

$$|\alpha_1-\alpha_2|^{(1-\tilde{\alpha}-\varepsilon)p} \leq \inf\left(\frac{c2^{(1-\tilde{\alpha}-\varepsilon)pt}}{\sup_{0\leq j\leq \tilde{t}} 2^{\frac{dpj}{2}}|1-2^{-jd}|}, \frac{c2^{(1-\tilde{\alpha}-\varepsilon)pt}}{\sup_{\tilde{\tilde{t}}<j\leq \tilde{t}} 2^{j(\frac{dp}{2}-d)}|2^{rt}2^{-jr(s+\frac{d}{2}-\frac{d}{r})}-1|}\right) \quad (55)$$

We have thus two cases:

- If $r > \frac{dp}{2s+d}$
$$\tilde{\alpha} = \frac{2s}{2s+d}.$$

But, if we take the first term, which satisfies
$$\sup_{0\leq j\leq \tilde{t}} 2^{\frac{dpj}{2}}|1-2^{-jd}| \sim 2^{\frac{tdp}{2s+d}},$$

we have
$$|\alpha_1 - \alpha_2|^{(1-\tilde{\alpha}-\varepsilon)p} \leq c2^{-\varepsilon pt}. \quad (56)$$

- When $r \leq \frac{dp}{2s+d}$, thus as $s > \frac{d}{r}$ we have necessarily $p > 2$ and
$$\tilde{\alpha} = \frac{2(s-\frac{d}{r}+\frac{d}{p})}{2(s-\frac{d}{r})+d}. \quad (57)$$

In this case,
$$\sup_{\tilde{\tilde{t}}<j\leq \tilde{t}} 2^{j(\frac{dp}{2}-d)}|2^{rt}2^{-jr(s+\frac{d}{2}-\frac{d}{r})}-1| \sim 2^{\frac{td(p-2)}{2(s-\frac{d}{r})+d}}.$$

And once again,
$$\forall t > 0 \quad |\alpha_1 - \alpha_2|^{(1-\tilde{\alpha}-\varepsilon)p} \leq c2^{-\varepsilon pt}. \quad (58)$$

As $1 - \tilde{\alpha} - \varepsilon > 0$, it can be deduced from equations (56) and (58) that for $t$ large enough, $M$ is of vanishing Lebesgue measure and $W((1-\tilde{\alpha}-\varepsilon)p,p)$ is an Haar null set in $B_r^{s,\infty}$. □

Thanks to invariance under inclusion, we have obtained that for every $\varepsilon > 0$, the set of functions $f$ in $B_r^{s,\infty}$ such that

$$\exists c > 0 \; \forall n \in \mathbb{N}, \; \mathbb{E}(\|\hat{f}_n^T - f\|_{L^p}^p) < c\sqrt{\frac{n}{\log n}}^{-(\alpha(s)+\varepsilon)p} \quad (59)$$

is a Haar null set.

The end of the proof follows similarly as for linear estimators.

## 4.4 Shrinkage procedures

The result of Theorem 2 is straightforward from Proposition 6 together with Propositions 7 and 8.